\documentclass[11pt]{article}
\usepackage{amsmath,amsfonts, amssymb,amsthm,multicol,graphicx,tikz, enumitem}
\usepackage{makeidx}
\renewcommand{\a}{\alpha}
\renewcommand{\b}{\beta}

\renewcommand{\l}{\lambda}

\newcommand{\ep}{\epsilon}
\newcommand{\ZZ}{\mathbb{Z}}
\newcommand{\RR}{\mathbb{R}} 
\newcommand{\PP}{\mathbb{P}}

\newcommand{\cE}{\mathcal{E}}
\newcommand{\cL}{\mathcal{L}}

\newcommand{\EE}{\mathbb{E}}

\newtheorem{prop}{Proposition}   
\newtheorem{lem}[prop]{Lemma} 
\newtheorem{thm}[prop]{Theorem} 
\newtheorem{corollary}[prop]{Corollary}

\title{A Critical Probability for Biclique Partition of $G_{n,p}$}
\author{Tom Bohman\thanks{This work was supported by a grant from the Simons
Foundation (587088, TB)} \, and Jakob Hofstad}

\begin{document}

\maketitle

\begin{center}
    {\bf Abstract}
\end{center}

The biclique partition number of a graph $G= (V,E)$, denoted $bp(G)$, is the minimum number of pairwise edge disjoint complete bipartite subgraphs of $G$ so that each edge of $G$ belongs to exactly one of them. It is easy to see that $ bp(G) \leq n - \a(G)$, where $\a(G)$ is the maximum size of an independent set of $G$. Erd\H{o}s conjectured in the 80's that for almost every graph $G$ equality holds; i.e., if $ G=G_{n,1/2}$ then $bp(G) = n - \a(G)$ with high probability. Alon showed that this is false. We show that the conjecture of Erd\H{o}s {\it is} true if we instead take $ G=G_{n,p}$, where $p$ is constant and less than a certain threshold value $p_0 \approx 0.312$. This verifies a conjecture of Chung and Peng for these values of $p$. We also show that if $p_0 < p <1/2$ then $bp(G_{n,p}) = n - (1 + \Theta(1)) \a(G_{n,p})$ with high probability.

\section{Introduction}
The biclique partition number of a graph $G=(V,E)$, denoted $bp(G)$, is the minimum number of pairwise edge disjoint complete bipartite subgraphs of $G$ so that each edge of $G$ belongs to exactly one of them. A well-known theorem of Graham and Pollak \cite{G-P} asserts that $bp(K_n) = n-1$, see \cite{Tverberg}, \cite{peck1984new}, \cite{Vishwanath} for more proofs.

Note that $bp(G) \le n - \a(G)$, where $\a(G)$ is the maximum size of an independent set of $G$. Indeed, one can partition all edges of $G$ into $n - \a(G)$ stars centered at the vertices of the complement of a maximum independent set of $G$. Erd\H{o}s  conjectured (see \cite{conjecture?}) that for almost every graph $G$, equality holds; i.e., that for the random graph $G = G_{n,1/2}$, $bp(G) = n - \a(G)$ with high probability (whp), that is, with probability that tends to 1 as $n$ tends to infinity. Chung and Peng \cite{chungpeng} extended the conjecture for the random graphs $G_{n, p}$ with $p \leq 0.5$, conjecturing that for any constant $p \leq 1/2$ we have $ bp(G) = n-(1+o(1))\a(G)$ whp. Alon proved in \cite{Alonbp} that Erd\H{o}s' conjecture for $ G_{n,1/2}$ is (slightly) incorrect, showing that for most values of $n$ we have $bp(G) \leq n - \a(G) - 1$ whp. A subsequent result of Alon, Bohman, and Huang \cite{ABHmore} improved the upper bound for $p = 0.5$: for some constant $c>0$, $bp(G) \leq n - (1+c) \a(G)$ whp, contradicting Chung and Peng's conjecture; however, they were not able to extend the result to $p < 1/2$.

Some progress has also been made on lower bounds for $bp(G_{n,p})$. Chung and Peng \cite{chungpeng} proved that, 
for $\Omega(1) \le p \le 1/2$ we have 
$$bp(G_{n,p}) \geq n - (\log n)^{3+\ep} = n -  (\log n)^{2 + \ep} \alpha( G_{n,p}). $$ 
Around the same time, Alon \cite{Alonbp} proved that if $\frac{2}{n} \leq p \leq c$ then
\begin{align*}
    bp(G_{n,p}) \ge n - O\left(\frac{\log(np)}{p}\right) = n - O \left( \alpha( G_{n,p}) \right),
\end{align*}
where $c>0$ is a small constant.

In this work we prove the theorem below, which defines a critical probability $p_0$ at which $ bp(G_{n,p})$ transitions from being equal to $ n - \a (G_{n,p})$ whp to being significantly less than this value whp. 
\begin{thm} \label{main}
Let $p_0$ be the real root of the polynomial $4x^3 - 7x^2 + 5x - 1$, so $p_0 \approx 0.312$. Then:
\begin{itemize}
    \item[i.] If $p \in (0,p_0)$ is a constant or $p = n^{\gamma}$ where $\gamma \in (-1/3,0)$ is a constant, then $bp(G_{n,p}) = n - \a(G_{n,p})$ whp. 
    \item[ii.] For any constant $p$ such that $p_0 \leq p < 1/2$, $bp(G_{n,p}) = n - O(\a(G_{n,p}))$ whp.
    \item[iii.] For any constant $p$ such that $p_0 < p < 1/2$, there exists some constant $c_p>0$ 
    such that $bp(G_{n,p}) \leq n - (1 + c_p)\a(G_{n,p})$ whp.
    \end{itemize}
\end{thm}
\noindent
Note that, if one combines parts ii and iii from the theorem for $ p_0 < p < 1/2$, then we have
\begin{align*}
    bp(G_{n,p}) = n - (1 + \Theta(1)) \a(G_{n,p}) \ \  whp. 
\end{align*}
Theorem~\ref{main} does not address $ bp( G_{n,1/2})$. Indeed, the question of whether or
not $ bp (G_{n,1/2}) = n - \Theta( \log n)$ whp remains an interesting open question. Furthermore,
Theorem~\ref{main} does not pin down $ bp( G_{n,p_0})$, and so the question of whether not
$ bp( G_{n,p_0})  = n - \alpha( G_{n,p_0}) $ whp is also interesting and open. These questions are discussed in
the Conclusion.

The proof of Theorem \ref{main} relies on proving the existence or non-existence of structures we call {\em special subgraphs of order $k$}. The existence of a special subgraph $G$ of order $k$ is equivalent to  $bp(G) \leq n-k$. The proof of part i uses the first moment method to bound the probability of existence of a special subgraph of order $\a(G_{n,p}) + 1$ and the fact that $ \a(G_{n,p})$ itself is concentrated on at most two values. (For constant $p$, two point concentration of the independence number goes back to early results of Bollob\'as and Erd\H{o}s \cite{1976cliques} and
Matula \cite{matula}.) The proof of part ii also uses the first moment method, but here we require a more sophisticated argument that 
also uses the BKR inequality regarding disjoint occurrence of events. The proof of part iii uses the second moment method to prove the existence of a higher-order special subgraph with some additional structure, which is added so as to make the second moment calculations easier; this closely follows methods introduced in \cite{ABHmore}. 

The remainder of the paper is organized as follows.  In the next section we establish some preliminaries, including the definition and discussion of special subgraphs and the BKR inequality. Section~2 also includes an explanation of the location of the critical value and a technical lemma that is crucial for the proofs of both part~i and part~ii of Theorem~\ref{main}. The proofs of parts i, ii, and iii of Theorem~\ref{main} appear in Sections 3, 4, and 5, respectively. We conclude the paper with some discussion of directions for further research in Section 6.  Unless explicitly stated otherwise we assume $ p < 1/2$ is a constant and $ G = G_{n,p}$.

\section{Preliminaries}
In following subsections we introduce the 
definition of a special subgraph,  explain the location of the critical probability $ p_0$, review the BKR inequality, and state and prove a technical lemma that
is crucial to parts i and ii of Theorem~\ref{main}.

\subsection{Special subgraphs}
We start with the definition of a structure that plays a central role throughout this paper.
For any graph $G$, let an induced subgraph $H$ of $G$ be 
called a {\em special subgraph of order $k$} if, for some $r \geq 0$, the
subgraph $H$ has $k + r$ vertices and the edge set of $H$ can be decomposed into at most $r$ non-star bicliques. 
Note that a special subgraph of order $k$ can simply be an independent set of size $k$ as we can have $r = 0$. The following lemma, which is analogous to Lemma 3.4 of \cite{Alonbp} and Lemma 12 of \cite{chungpeng}, establishes a connection between the existence of special subgraphs and $bp(G)$.
\begin{lem}
\label{special}
Let $2 \le k \le n$ be integers.
If $G$ is a graph on $n$ vertices, then $bp(G) \leq n-k$ if and only if $G$ has a special subgraph of order $k$.
\end{lem}
\begin{proof}
If $G$ has a special subgraph $H$ of order $k$, then consider the biclique decomposition of $G$ given by the decomposition of $H$ into $r$ bicliques together with $n-k-r$ stars centered at the other vertices. Now assume that $bp(G) \leq n-k$. Consider a partition of $G$ into at most $n-k$ bicliques. Iteratively delete the central vertices of star components until no star components remain (Note that this may involve deleting edges from non-star bicliques that later become stars). If $s$ vertices are deleted in this process, then the induced subgraph that remains is a graph $G'$ on $n-s$ vertices whose edge set is partitioned into at most $ n-k-s$ bicliques, each of which is not a star. 
It follows that $G'$ is a special subgraph of order $k$, where $r = n-k-s$. 
\end{proof}

\subsection{The critical value}
Next, we state and discuss a lemma that highlights the technical role of the value of $p_0$. The proof of this
fact is left as a calculus exercise.
\begin{lem} \label{root}
Let $p_0$ be the real root of the polynomial $4x^3 - 7x^2 + 5x - 1$. If $p \in (0,1)$, then we have
\begin{align*}
     1+\left(\frac{p}{1-p}\right)^2 - (1 - p)^{-1/2} \ \ \begin{cases} < 0 &  \text{ if } p < p_0 \\
     > 0 & \text{ if } p > p_0.
     \end{cases}
\end{align*}
\end{lem}

In order to explain the role of the expression in Lemma~\ref{root} in the arguments below, 
we need to make a definition. We say that a biclique defined by vertex parts $ B_1,B_2 $ is {\em based} on vertex set $B_i$ if 
$ |B_i| = \min \{ |B_1|, |B_2| \}$. (N.b. If $|B_1| = |B_2|$ then the biclique is based on two sets.) 
It turns out that bicliques that are based on sets of size two play the central role in the appearance of special subgraphs of order $k$ with $ k > \alpha( G_{n,p})$.  In particular, it turns out that a special subgraph of order $k$ with $ k > \alpha ( G_{n,p})$ appear whp if and only if 
a special subgraph of order $k$ consisting of bicliques based on sets of size 2 appear whp.
This is reflected in Lemmas~\ref{combo}~and~\ref{prob bound} below. In the remainder of this subsection we give a rough account of how this key observation regarding bicliques based on sets of size two yields the critical value defined in Lemma~\ref{root}.

Assume for the sake of this discussion that we only need to
consider bicliques based on sets of size 2. We first note
that the probability that a vertex set $K$ induces a special subgraph that contains
a single biclique based on a vertex set $B$ of size 2 is roughly 
\[ \sum_{B' \subset K \setminus B} p^{ 2 |B'|} (1-p)^{ \binom{ |K|}{2} - 2 |B'| } 
= (1-p)^{ \binom{|K|}{2}}  \left[ 1 + \left( \frac{p}{ 1-p} \right)^2 \right]^{|K|-2}.
\]
Extending this observation to multiple bicliques based on sets of size 2, we see that 
the expected number of special subgraphs of order $k$ 
can be roughly approximated by
\begin{multline*}
 \sum_{ r \ge 0} \binom{n}{ k+r} (1-p)^{ \binom{k+r}{2}} \left( \binom{k+r}{2} \left[ 1 + 
 \left( \frac{p}{ 1-p} \right)^2 \right]^{k+r-2} \right)^r   \\
 \approx \binom{n}{k} (1 -p)^{\binom{k}{2}} \sum_{r\ge 0} \left( n (1-p)^{k +r/2}  \left[ 1 + 
 \left( \frac{p}{ 1-p} \right)^2 \right]^{k+r} \right)^r.
\end{multline*}
(Note that the first step in this estimate makes crucial use of the fact that the bicliques
are edge-disjoint.)
Now, recall that, restricting our attention to constant $p$, we have 
$ \alpha( G_{n,p}) = (2+o(1))\log_b(n)$, where $ b = \frac{1}{1-p}$, 
and if $k$ is larger than $ \alpha( G_{n,p})$ 
then the expression before the summation, which is the expected number of independent
sets of size $k$, is small. Furthermore, as $ (1-p)^k = n^{-2 +o (1)}$, whether
or not the expected number of special subgraphs of order $k$ vanishes depends on the
relationship between $ 1 + \left( \frac{p}{ 1-p} \right)^2$ and $ ( 1-p)^{-1/2}$. This is why
the critical value is at the point $p_0$ defined in Lemma~\ref{root} above. Of course, these rough
calculations need to be developed into a rigorous proof. This is done in the sections below; in particular,
calculations along these lines appear in the arguments leading up to (\ref{final_W_bound}) and (\ref{plug in k}).

\subsection{The BKR inequality}
\label{bkr}
Consider a product probability space defined on $ \Omega = S^d$, where $S$ is a finite set.  
The BKR inequality deals with `disjoint' occurrences of events in this setting,
asserting that the probability of such a disjoint occurrence of events is bounded above by the product
of the probabilities of the events in question. This inequality is named in honor of van den Berg and
Kesten, who conjectured the inequality, and Reimer, who proved it.

In order to make this notion precise, we need to introduce some definitions. 
For $K \subseteq \{1, \dots, d\}$ and $ \omega = ( \omega_1, \dots, \omega_d) \in S^d$
we define the {\em cylinder} determined by $K$ at $ \omega$ as
\[ Cyl(\omega,K) = \{ (\alpha_1, \dots, \alpha_d) \in S^d: \alpha_i = \omega_i \ \forall i \in K \}. \]
We define the complement of the index set $K$ to be $K^c = \{1, \dots, d\} \setminus K$.  We say $ \omega$ is a
{\em disjoint occurrence} of events $A, B \subseteq S^d$ if there is some index set $K \subset \{1, \dots, d\}$ such that
$ Cyl(\omega,K) \subseteq A$ and $ Cyl(\omega, K^c) \subseteq B$.  In words, we say that $ \omega$ is a disjoint occurrence of 
$A$ and $B$ if $ \omega$ is in the intersection of these two events with the additional property that there is some
set of coordinates $K$ such that any occurrence that agrees with $ \omega$ on $K$ is in $A$ and any occurrence that agrees
with $ \omega$ on the complement of $K$ is in $B$. We use the notation $ A \square B$ to denote the collection of
disjoint occurrences of events $A$ and $B$. So, if $ A, B \subseteq S^d$, then we set
\[ A \square B = \{ \omega \in S^d: \ \exists K \subseteq \{1, \dots, d\} \text{ such that }  Cyl(\omega,K) \subseteq A \text{ and } Cyl(\omega,K^c) \subseteq B\}.   \]
\begin{thm}[BKR inequality. See \cite{preCylinder1,Cylinder}] \label{BKR}
If  $ A, B \subseteq S^d$, then
\[ \PP( A \square B) \le \PP(A) \PP(B). \]
\end{thm}
\noindent
The BKR inequality can be
extended to collections of events with a suitable definition of disjoint
occurrence of events (we follow the treatment of \cite{generalBKR2, generalBKR1}). 
We begin by noting that, for each event $A \subseteq S^d$ and index $ J \subseteq \{1, \dots, d\}$, there is some set of elements that gain entry into
$A$ by virtue of their values on the set $J$. This is the set
\[ [A]_J = \{ \omega \in A: Cyl(\omega,J) \subseteq A \}, \]
which may be empty. Note that
\[ A \square B = \bigcup_{ \stackrel{ J,K \subseteq \{1,\dots, d\} }{ \text{disjoint}}} [A]_J \cap [B]_K. \]
Now, for events $ A_1, \dots, A_r \subseteq S^d$, we define
\[ \square_{i=1}^r A_i = \bigcup_{ \stackrel{ J_1 ,\dots ,J_r \subseteq \{1,\dots, d\} }{ \text{disjoint}}}  \bigcap_{i=1}^r [A_i]_{J_i}. \]
Note that a disjoint occurrence of the events $ A_1,\dots, A_r$ is an element of the intersection of these events
with the property that there is some collection of $r$ disjoint subsets $ J_1, \dots, J_r$ of the index set with the property that
$\omega$ gains entry into each $A_i$ by virtue of its values on the corresponding set $J_i$. Given this definition, the BKR
inequality generalizes as follows:
\begin{thm} \label{BKRmany}
If  $ A_1, \dots, A_r \subseteq S^d$, then
\[ \PP\left( \square_{i=1}^r A_i \right) \le  \prod_{i=1}^r \PP(A_i). \]
\end{thm}
\noindent
Theorem~\ref{BKRmany} follows from Theorem~\ref{BKR} by a simple induction argument.
For details and further discussion see \cite{generalBKR2}.

\subsection{Two technical lemmas}
We conclude this section with two technical lemmas that are crucial to both parts i and ii of Theorem~\ref{main}. For integers $ r > s \ge 0$ define
$$R(r,s,p) := \frac{r-s}{r+s+2} \left(\frac{1-p}{p}\right)^{s+1}.$$

\begin{lem} \label{lem: ratio large}
    If $p$ be a constant in $(0,1/2)$ then there is an integer $N_p$, dependent only on $p$, such that, for any integer $r \geq N_p $ and any $s \in \{0, 1, 2, \dots, r-1\}$, we have $R(r,s,p) > 1$. Furthermore, if $ p < 0.01$ then $R(r,s,p) > 2$.
\end{lem}
\begin{proof}
    First, we show that $R(r,s,p) > 1$ for constant $p$ and large enough $r$. Let $ \xi = (1-p)/p$. Note that $\xi$ is a constant and $ \xi > 1$. We have two cases:
    \begin{enumerate}
        \item If $s \leq \sqrt{r}$ then we have $\frac{r-s}{r+s+2} \geq 1-\frac{2\sqrt{r}+2}{r+\sqrt{r}+2}$, which tends to 1 as $r \to \infty$. Note that $ \xi^{s+1} \geq \xi > 1$, which is a fixed constant. Thus, we have 
        $R(r,s,p) > 1$ for large enough $r$.
        \item For $s > \sqrt{r}$ then we have $\frac{r-s}{r+s+2} \geq \frac{1}{2r+1}$ and 
        $\xi^{s+1} > \xi^{\sqrt{r}}$. As $ \xi^{s+1}$ is at least exponential in $\sqrt{r}$ 
        we have $R(r,s,p) > 1$ for $r$ sufficiently large. 
    \end{enumerate}
    Next, we take away the restriction $r \geq N_p$ and assume $p < 0.01$. We again consider two cases:
    \begin{enumerate}
        \item If $s \leq \frac{r}{2}$, then $\frac{r-s}{r+s+2} \geq \frac{1}{5}$, and $\xi^{s+1} \geq \xi > 10$. Thus, $ R(r,s,p) > 2$.
        \item If $s > \frac{r}{2}$, then $\frac{r-s}{r+s+2} > \frac{1}{3r} \geq \frac{1}{3^r}$ while $\xi^{s+1} > \xi^{r/2}$. Hence we have $R(r,s,p) > \left(\xi/9\right)^{r/2} > 2$ for all $r \in \ZZ^+$.
    \end{enumerate}
\end{proof}

\begin{lem} \label{combo} 
If $m' \geq 2$ is an integer and $ p \in (0,1/2)$ is a constant then 
there exists an integer $M_{p,m'}$, dependent only on $p$ and $m'$, such that if $n' > M_{p,m'}$ is an integer then we have
\begin{align} \label{crux}
    \sum_{\substack{ \l \geq \mu \geq m' \\ \l + \mu \leq n'}} \binom{n'}{\l+\mu} \binom{\l+\mu}{\l} \left(\frac{p}{1-p} \right)^{\l \mu} \leq 2 \binom{n'}{m'}\left(1 + \left(\frac{p}{1-p}\right)^{m'} \right)^{n'
    }.
\end{align}
 Furthermore, if $ p < 0.01$ then inequality (\ref{crux}) holds for all integers $m' \geq 2$  and $n' \ge 2m'$.
\end{lem}

\begin{proof} First, we use the Binomial Theorem to bound the sum if we fix $\mu = m'$:
\begin{align}
    \sum_{\l = m'}^{n' - m'} \binom{n'}{\l+m'} \binom{\l+m'}{\l} \left(\frac{p}{1-p} \right)^{\l m'} &= \binom{n'}{m'} \sum_{\l = m'}^{n' - m'} \binom{n'-m'}{\l} \left(\frac{p}{1-p} \right)^{\l m'} \nonumber \\&\leq 
    \binom{n'}{m'}\left(1 + \left(\frac{p}{1-p}\right)^{m'} \right)^{n'}. \label{equ: binom thm}
\end{align}
Next, let $r = \l + \mu$ and let $s = \l - \mu$. Then $\l = \frac{r + s}{2}$ and $\mu = \frac{r - s}{2}$. We now rewrite the sum in terms of indices $r$ and $s$:
\begin{multline}
    \sum_{ \substack{ \l \geq \mu \geq m' \\ \l + \mu \leq n'}} \binom{n'}{\l+\mu} \binom{\l+\mu}{\l} \left(\frac{p}{1-p} \right)^{\l \mu} \\=
    \sum_{r = 2m'}^{n'} \ \  \sum_{ \substack{ 0 \leq s \leq r-2m'\\ r + s \text{ is even}} } \binom{n'}{r} \binom{r}{(r+s)/2}\left(\frac{p}{1-p} \right)^{(r/2)^2 - (s/2)^2}. \label{the rest prime}
\end{multline}
Next, let
\begin{align*}
    f_{n',r,s} := \binom{n'}{r} \binom{r}{(r+s)/2}\left(\frac{p}{1-p} \right)^{(r/2)^2 - (s/2)^2}.
\end{align*}
We want to sum $f_{n',r,s}$ over all feasible $r$ and $s$. Observe that $\frac{f_{n',r,s+2}}{f_{n',r,s}} = R(r,s,p)$, where $R(r,s,p)$ is defined before Lemma~\ref{lem: ratio large}. \\

We now take the two separate given conditions in the Lemma one by one: first, assume $p < 0.01$. 
By~Lemma \ref{lem: ratio large}, if $s < r$ then $R(r,s,p) > 2$, therefore
\begin{align}
    \sum_{r = 2m'}^{n'} \ \  \sum_{ \substack{ 0 \leq s \leq r-2m'\\ r + s \text{ is even}} } f_{n',r,s} &\leq  2 \sum_{r = 2m'}^{n'} f_{n',r,r-2m'} \label{equ: first one} \\&=
     2\sum_{\l = m'}^{n' - m'} \binom{n'}{\l+m'} \binom{\l+m'}{\l} \left(\frac{p}{1-p} \right)^{\l m'} \nonumber \\&\leq
     2 \binom{n'}{m'}\left(1 + \left(\frac{p}{1-p}\right)^{m'} \right)^{n'
    } \label{equ: second one} \qquad \text{by (\ref{equ: binom thm}).}
\end{align}

Next, consider $p \in (0,1/2)$ fixed and let $n' > M_{p,m'}$ for $M_{p,m'}$ to be defined later; in other words, take $n'$ to be ``arbitrarily large" depending on just $m'$ and $p$. Note that the values of $f_{n',r,s}$ for the first few values of $r$ will be asymptotically smaller (in terms of $n'$) than those values of $f_{n',r,s}$ for higher values of $r$, since $\binom{n'}{r}$ is the dominant factor; thus, for large enough $n'$ (with $N_p$ also satisfying Lemma \ref{lem: ratio large}):
\begin{equation*}
    \sum_{r = 2m'}^{n'} \ \sum_{ \substack{ 0 \leq s \leq r-2m'\\  r + s \text{ is even}}} f_{n',r,s} \leq \frac32 \sum_{r = N_p}^{n'} \ \sum_{ \substack{ 0 \leq s \leq r-2m' \\  r + s \text{ is even}}} f_{n',r,s}.
\end{equation*}
By Lemma \ref{lem: ratio large}, we have $R(r,s,p) > 1$ for $r \geq N_p$, hence
\begin{align*}
& \frac32 \sum_{r = N_p}^{n'} \ \  \sum_{ \substack{ 0 \leq s \leq r-2m'\\ r + s \text{ is even}} } f_{n',r,s} \\&= \frac32\sum_{r = N_p}^{n'} f_{n',r,r-2m'} + \frac32\sum_{r = N_p}^{n'} \ \  \sum_{ \substack{ 0 \leq s \leq r-2m'-2\\ r + s \text{ is even}} } f_{n',r,s} \\&<
\frac32\sum_{r = N_p}^{n'} f_{n',r,r-2m'} + \sum_{r = N_p}^{n'} n' f_{n',r,r-2m'-2} \\&\leq
\frac32\binom{n'}{m'}\left(1 + \left(\frac{p}{1-p}\right)^{m'} \right)^{n'} + n' \binom{n'}{m'+1}\left(1 + \left(\frac{p}{1-p}\right)^{m'+1} \right)^{n'},
\end{align*}
where the last inequality follows from (\ref{equ: binom thm}) as in (\ref{equ: first one})-(\ref{equ: second one}), replacing $m'$ with $m'+1$ to get the second term. Since $p <1/2$, 
the first summand is asymptotically larger than the second summand as $n' \to \infty$ (because the exponential factor has a larger base and recalling that $m'$ is fixed).  Therefore for some integer $M_{p,m'}$ we have
\begin{equation*}
    \sum_{ \substack{ \l \geq \mu \geq m' \\ \l + \mu \leq n'}} \binom{n}{\l+\mu} \binom{\l+\mu}{\l} \left(\frac{p}{1-p} \right)^{\l\mu} 
    < 2 \binom{n'}{m'}\left(1 + \left(\frac{p}{1-p}\right)^{m'} \right)^{n'}
\end{equation*}
for all $n' > M_{p,m'}$,
as desired.
\end{proof}

\section{Lower bound on $ bp( G_{n,p})$ for $p<p_0$}

In this section we prove part i of Theorem~\ref{main}.
Applying Lemma~\ref{special}, our objective is to show that special subgraphs of order $ k = \a +1$ do not exist whp.

We begin with a lemma that bounds the probability that a
fixed set of vertices in the random graph gives a special subgraph. This lemma is crucial in the first moment calculations below. Recall that $ M_{p,2}$ is defined in Lemma~\ref{combo}.
\begin{lem} \label{prob bound}
If $n' > M_{p,2}$, the probability that $E(G_{n',p})$ can be decomposed into at most $r$ nonempty non-star bicliques is at most
\begin{align*}
    2 (1-p)^{\binom{n'}{2}} \left[2 \binom{n'}{2}\left(1 + \left(\frac{p}{1-p} \right)^2 \right)^{n'} \right]^r.
\end{align*}
If $ p < 0.01 $ then the same bound holds for all $n' \geq 4$.
\end{lem}
\begin{proof}
We begin with the case that $E(G_{n',p})$ can be decomposed into exactly $r$ nonempty non-star bicliques. Let $P_r$ be the probability of this event.
Consider a fixed graph $G$ on $n'$ vertices with such a decomposition, and let $H_1,\dots,H_r$ be the bipartite graphs in the decomposition. For each $H_i$, let $A_i$ and $B_i$ be the parts of the bipartition, where $|A_i| = \l_i$ and $|B_i| = \mu_i$, and assume wlog 
that $\l_i \geq \mu_i$. Then $|E(H_i)| = \l_i \mu_i$. We also have $\l_i, \mu_i \geq 2$ and $\l_i + \mu_i \leq n'$.
Now bring in $G_{n',p}$. We have that, for {\it labeled} graphs $G_{n',p}$ and $G$:
\begin{align*}
    \PP[G_{n',p} = G] = (1 - p)^{\binom{n'}{2}} \left(\frac{p}{1-p}\right)^{\l_1
    \mu_1 + \l_2\mu_2 + \dots + \l_r \mu_r},
\end{align*}
since we need to choose $\l_1 \mu_1 + \l_2 \mu_2 + \dots + \l_r \mu_r$ specific edges and have all other pairs of vertices not be edges. Now consider the number of choices of parts $A_i$ and $B_i$ if we fix each $\l_i$ and $\mu_i$. For a single bipartite graph $H_i$, we have $\binom{n'}{\l_i + \mu_i} \binom{\l_i + \mu_i}{\l_i}$ ways to choose its partite sets (this is a double count if $\l_i = \mu_i$). Thus, an upper bound for the number of choices of all the bipartite graphs will be $\prod_{i=1}^{r}\binom{n'}{\l_i + \mu_i} \binom{\l_i + \mu_i}{\l_i}$. Therefore, we have
\begin{align*}
    P_r &\leq \sum_{ \substack{ \l_i \geq \mu_i \geq 2\\ \l_i + \mu_i \leq n'}} \left(\prod_{i=1}^{r}\binom{n'}{\l_i + \mu_i} \binom{\l_i + \mu_i}{\l_i} \right)(1 - p)^{\binom{n'}{2}} \left(\frac{p}{1-p}\right)^{\l_1\mu_1 + \l_2\mu_2 + \dots + \l_r \mu_r} \\&=
    (1 - p)^{\binom{n'}{2}} \sum_{ \substack{ \l_i \geq \mu_i \geq 2 \\ \l_i + \mu_i \leq n'}} \prod_{i=1}^{r}\binom{n'}{\l_i + \mu_i} \binom{\l_i + \mu_i}{\l_i} \left(\frac{p}{1-p}\right)^{\l_i \mu_i} \\&=
    (1 - p)^{\binom{n'}{2}} \left[ \sum_{ \substack{ \l \geq \mu \geq 2 \\ \l + \mu \leq n'}} \binom{n'}{\l + \mu} \binom{\l + \mu}{\l} \left(\frac{p}{1-p}\right)^{\l \mu}\right]^r
    \\& \le (1-p)^{\binom{n'}{2}} \left[2 \binom{n'}{2}\left(1 + \left(\frac{p}{1-p} \right)^2 \right)^{n'} \right]^r,
\end{align*} 
noting that we apply Lemma~\ref{combo} with  $m'=2$.

Now consider the case that $ G$ is decomposed into at most $r$ bicliques. We consider the sum of the above bounds on $ P_{r'}$ over $r'$ between 0 and $r$. This is just the sum of a geometric sequence with a ratio greater than 2, and hence bounded above by twice the value of the last summand.
\end{proof}

We now have the tools to prove part i of Theorem \ref{main}. We
break the argument into two parts depending on whether or not $p$ goes to zero with $n$. For $ p \in (0,p_0)$ a constant, our argument follows the structure of the classical argument that establishes two point concentration of the independence number for $p$ in this range (see Section 11.1 of \cite{bela} for details).

\begin{prop} \label{prop_constant_p}
Let $p \in (0,p_0)$ be constant.
There exists an integer function $k = k(n,p)$ such that the following three events happen with high probability:
\begin{enumerate}[label=\roman*.]
    \item There exists an independent set of size $k-1$. 
    \item There does not exist any special subgraph of order $k$ that is not an independent set.
    \item There does not exist any special subgraph of order $k+1$.
\end{enumerate}
\end{prop}
\noindent
The fact that the independence number of $ G_{n,p}$ is concentrated on the values $k, k-1$ (where $k$ is 
defined in the proof below) is a well known. See, for example, the original proofs of Bollob\'as and Erd\H{o}s \cite{1976cliques} or Matula \cite{matula} or Theorem~4.5.1 in \cite{probmethod}. We include a proof of part i of Proposition~\ref{prop_constant_p} for completeness and to set the stage for the corresponding result for $p = n^{\gamma}$ with $ -1/3 < \gamma < 0$.
\begin{proof}
First, we define $k$ by first defining variables $b$ and $\ep_p$. Set $b = 1/(1-p)$, and let $\ep_p \in (0,p)$ be a constant such that
\begin{align}
    1 + \left(\frac{p}{1-p}\right)^2 < (1-\ep_p)(1-p)^{-1/2} \label{ep_p_for_constant_p},
\end{align}
which is possible by Lemma \ref{root}. Next, let $k$ be the smallest positive integer such that
\begin{align}
    \binom{n}{k} (1-p)^{\binom{k}{2}} < n^{-\log_b(1-\ep_p)}.  \label{k_p_constant}
\end{align}
Note that $n^{-\log_b(1-\ep_p)} = \omega(1)$, and that \vspace{-3mm}

\begin{equation}
    k = (2 -o(1)) \log_b(n). \label{equ: k first approx}
\end{equation}

We use the second-moment method to prove (i). Note that it follows from (\ref{k_p_constant}) that the expected number of independent sets of size $k-1$ is $ \omega(1)$. For any induced subgraph $K$ of $G$ with $|V(K)| = k-1$, let $X_K$ be the indicator random variable of the event that $K$ is an independent set, and let $X = \sum_{K} X_K$ be the random variable that stands for the number of independent sets of size $k-1$ in $G$. Then $\EE(X) = \omega(1)$. By Chebyshev's Inequality, it is enough to show that Var$(X)/\EE(X)^2 = o(1)$. 

For $K,K' \subset V(G),|K|=|K'|=k-1$, let $K \sim K'$ denote that $2 \leq | K \cap K'| \leq k-2$. The variance of $X$ satisfies
\begin{align*}
    \text{Var}(X) = \sum_K \text{Var}(X_K)+ \sum_{K \sim K'} \text{Cov}(X_K,X_{K'}) \leq \EE(X)+ \sum_{K \sim K'} \EE(X_K X_{K'}).
\end{align*}
For each $i \in [2, k - 2]$, let $h_i$ denote the contribution of the pairs with intersection size $i$ to the
above sum divided by $\EE(X)^2$, that is:
\begin{align*}
    h_i = \sum_{|K \cap K'| = i} \frac{\EE(X_K X_{K'})}{\EE(X)^2}.
\end{align*}
Our current objective is to show that $\sum_{i=2}^{k-2} h_i = o(1)$; this, combined with $\EE(X) = \omega(1)$, gives us (i).

It is straightforward to calculate $h_i$ for each $i \in [2,k-2]$:
\begin{align*}
    h_i = \frac{\binom{k-1}{i} \binom{n-k+1}{k-i-1}}{\binom{n}{k-1}} (1-p)^{-\binom{i}{2}}.
\end{align*}
Next, for $i \in [2,k-3]$, let $\kappa_i := h_{i+1}/h_i$. Then we have
\begin{align*}
    \kappa_i = \frac{(1-p)^{-i}(k-i-1)^2}{(i+1)(n -2k+i+3)} = n^{2(i/k)-1+o(1)}.
\end{align*}
It is easy to validate that the dominant terms in the sequence $h_2,h_3,\dots,h_{k-2}$ will be $h_2$ and $h_{k-2}$, so we have
\begin{align}
    \sum_{i=2}^{k-2} h_i &= (1+o(1)) (h_2 + h_{k-2}) \nonumber \\&=
    (1+o(1)) \left(\frac{\binom{k-1}{2} (k-1)(k-2) \prod_{j=1}^{k-3} (n - 2k + j +4 )}{(1-p)\prod_{j=1}^{k-1}( n - k + j + 1)} + \frac{(k-1)(n-k+1)}{\binom{n}{k-1}} (1-p)^{-\binom{k-2}{2}}\right) \nonumber \\&\leq
    (1+o(1)) \left(\frac{k^4}{2(1-p) n^2} + \frac{kn(1-p)^{k-2}}{\EE(X)}\right) \label{littleo} \\&\leq
    n^{-2+o(1)} + n^{1-2+o(1)}/\EE(X) \nonumber \\&=
    o(1). \nonumber
\end{align}
Thus, (i) holds w.h.p.

Using a straightforward first-moment calculation, we prove that both (ii) and (iii) hold whp. We will prove (ii), then explain briefly how (iii) quickly follows. First, let $W$ be the random variable defined as the number of nonempty special subgraphs of order $k$ in $G$. To prove that (ii) holds whp, it is enough to show that $\EE(W) = o(1)$ by Markov's Inequality. By Lemmas \ref{root} and \ref{prob bound}, for large enough $n$ we have
\begin{align}
    \EE(W) &\leq \sum_{r=1}^{n-k} \binom{n}{k+r} 2(1-p)^{\binom{k+r}{2}}\left[2 \binom{k+r}{2}\left(1 + \left(\frac{p}{1-p} \right)^2 \right)^{k+r} \right]^{r} \nonumber \\&\leq
    2\binom{n}{k} (1-p)^{\binom{k}{2}}\sum_{r=1}^{n-k} \left(\frac{n}{k}\right)^r (1-p)^{(k + r/2 - 1/2)r}\left[2 \binom{k+r}{2}\left(1 + \left(\frac{p}{1-p} \right)^2 \right)^{k+r} \right]^{r} \nonumber \\&=
    2\binom{n}{k} (1-p)^{\binom{k}{2}}\sum_{r=1}^{n-k}\left[2\left(\frac{n}{k}\right) (1-p)^{k + r/2 - 1/2} \binom{k+r}{2}\left(1 + \left(\frac{p}{1-p} \right)^2 \right)^{k+r} \right]^{r} \label{W bound} \\&\leq 
    2\binom{n}{k} (1-p)^{\binom{k}{2}}\sum_{r=1}^{n-k}\left[2 \left(\frac{n}{k}\right) (1-p)^{k/2 - 1/2}(1-\ep_p)^{k+r} \binom{k+r}{2} \right]^{r}. \nonumber
\end{align}
For large enough fixed $k$ and $n$, the expression inside the square brackets is maximized at $r=1$, and would be increased if $r$ was made equal to 0 (by the ratio test), so
\begin{align}
    \EE(W) \leq
    2\binom{n}{k} (1-p)^{\binom{k}{2}}\sum_{r=1}^{n-k}\left[ n (1-p)^{k/2 - 1/2}(1-\ep_p)^{k} k \right]^{r}. \label{final_W_bound}
\end{align}
By (\ref{k_p_constant}) and (\ref{equ: k first approx}), we have 
\begin{align*}
    \EE(W) &\leq n^{-\log_b(1-\ep_p)} \sum_{r=1}^{n-k}\left[n^{2\log_b(1-\ep_p)+o(1)}\right]^r \\&=
    o(1),
\end{align*}
thus (ii) holds whp.

We bound the probability of (iii) by letting $W'$ be the variable that counts the number of special subgraphs of order $k+1$ in $G$ (which are not necessarily nonempty). The expected value calculations are the same, with the exception that $k$ changes to $k+1$ and the summation starts at $r=0$ rather than $r=1$. 
Hence we can borrow (\ref{final_W_bound}) to bound $\EE[W']$:
\begin{align*}
    \EE(W') \leq
    2 \binom{n}{k+1} (1-p)^{\binom{k+1}{2}}\sum_{r=0}^{n-k-1}\left[ n (1-p)^{(k+1)/2 - 1/2}(1-\ep_p)^{k+1} (k+1) \right]^{r}.
\end{align*}
The main change in the calculations is the new `independence set'
part $\binom{n}{k+1}(1-p)^{\binom{k+1}{2}}$. Note that we have
\begin{align*}
    \binom{n}{k+1} (1-p)^{\binom{k+1}{2}} &= \binom{n}{k} (1-p)^{\binom{k}{2}} \left(\frac{(1-p)^{k} (n-k)}{k+1}\right) \\&\leq
    n^{-\log_b(1-\ep_p)} \left(\frac{ n^{-1+o(1)}}{k}\right) \\&=
    n^{-1-\log_b(1-\ep_p) + o(1)} \\&=
    o(1) \qquad  \text{ (since $\ep_p < p$)}.
\end{align*}
Since the expression inside the square brackets above is $n^{2 \log_b(1-\ep_p) + o(1)}$ like in (\ref{final_W_bound}), then $\EE(W') = o(1)$, so (iii) holds whp.
\end{proof}

\begin{prop} \label{prop_nonconstant_p}
Let $\gamma \in (-1/3,0)$ be a constant and set $p = n^{\gamma}$. Then there exists an integer function $k = k(n,\gamma)$ such that $ G_{n,p}$ has the following three properties with high probability:
\begin{enumerate}[label=\roman*.]
    \item There exists an independent set of size $k-1$.
    \item There does not exist any special subgraph of order $k$ that is not an independent set.
    \item There does not exist any special subgraph of order $k+1$.
\end{enumerate}
\end{prop}

\begin{proof}
This proof is very similar to the proof of Proposition \ref{prop_constant_p}.
First, we re-define $k$. Let $k$ be the smallest integer such that
\begin{align}
    \binom{n}{k} (1-p)^{\binom{k}{2}} < n^{\gamma + 1/3} \label{k_notconstant_p}.
\end{align}
Note that $ n^{\gamma+1/3}= \omega(1)$. Also, it is straightforward to show that 
\begin{align*}
k = \frac{(2 + o(1))(\ln(np))}{p} = \frac{(2 + o(1))(1+\gamma)\ln(n)}{p}.
\end{align*}

The proof for (i) is similar to the corresponding argument in the proof of Proposition~\ref{prop_constant_p}. We define $X$ and $h_2, \dots, h_{k-2}$ in the same way. Again, $\EE(X) = \omega(1)$, so again it suffices to show that $\sum_{i=2}^{k-2} h_i = o(1)$, where for each $i \in [2,k-2]$, we have
\begin{align*}
    h_i = \frac{\binom{k-1}{i} \binom{n-k+1}{k-i-1}}{\binom{n}{k-1}} (1-p)^{-\binom{i}{2}}.
\end{align*}
For any $ \epsilon >0$ and $ \epsilon k < i < (1-\epsilon) k$, we have
\begin{equation*}
    \begin{split}
    h_i & \le 2^k \left( \frac{k}{n} \right)^i (1-p)^{-\binom{i}{2}} \\
    & \le \left[  2^{k/i} n^{-1 - \gamma +o(1)} e^{pi/2}   \right]^i \\        
    & \le \left[  2^{ 1/\epsilon}  n^{-1 - \gamma +o(1)} \cdot n^{ (1 -\epsilon) ( 1 + \gamma) +o(1)} \right]^i\\
    & = \left[  n^{ -\epsilon( 1 + \gamma) + o(1) } \right]^i,
    \end{split}
\end{equation*}  
hence $\sum_{i=\ep k}^{(1-\ep) k} h_i = o(1)$.
For the remaining values of $i \in [2,k-3]$, consider $\kappa_i = h_{i+1}/h_i$. As before,
we have
\begin{align*}
    \kappa_i = \frac{(1-p)^{-i}(k-i-1)^2}{(i+1)(n -2k+i+3)}.
\end{align*}
For $ i < \epsilon k$ we have
\[ \kappa_i \le n^{2 \epsilon( 1 + \gamma) + o(1)} \cdot n^{-1 - 2\gamma} = o(1) \]
for $n$ sufficiently large, provided $ \epsilon < ( 1 + 2\gamma)/( 2 + 2\gamma) $.  
Thus, we can bound the sum of the first $ \epsilon k$
terms in $ \sum h_i$ by $ (1 + o(1)) h_2$ using a geometric series.  
Furthermore, if $ i > (1 -\epsilon)k$
then we have
\begin{equation}
\label{eq:adjust1}
\kappa_i \ge n^{ (1 -\epsilon )(2 + 2 \gamma) + o(1)} \cdot n^{-1 + \gamma}
= n^{ 1 - 2\epsilon + \gamma(3 - 2\epsilon)+ o(1) }= \omega(1), 
\end{equation}
for $ \epsilon $ sufficiently small. 
This allows us to bound the sum of the final $ \epsilon k$ terms by 
$ (1 + o(1)) h_{k-2}$.  Putting these observations together and borrowing $(\ref{littleo})$, we have
\begin{align}
    \sum_{i=2}^{k-2} h_i &= o(1) + (1+o(1)) (h_2 + h_{k-2}) \nonumber \\
    &\leq
    o(1)+(1+o(1)) \left(\frac{k^4}{2(1-p) n^2} + \frac{k(n-k)(1-p)^{k-2}}{\EE(X)}\right) \nonumber \\
    & \leq o(1)+n^{-2 - 4\gamma + o(1)} + n^{1- \gamma} \cdot n^{-2( 1 + \gamma) + o(1)} / n^{\gamma + 1/3} \label{eq:adjust2} \\
    & = o(1) + n^{ -4/3 - 4\gamma +o(1)} \nonumber \\
    & = o(1) \nonumber
\end{align}
as desired.

To prove (ii), we again let the variable $W$ be the number of nonempty special subgraphs of order $k$. As Lemma~\ref{prob bound} can be used for $p = o(1)$, we can borrow the calculations from the proof of Proposition \ref{prop_constant_p} up until (\ref{W bound}). In the expression that follows, we replace $\ep_p$ with $p/2 - p^2$, which satisfies (\ref{ep_p_for_constant_p}) for large enough $n$ (e.g. consider the power series of $\left(1 + (p/(1-p)\right)^2)(1-p)^{1/2}$). This leaves us with
\begin{align}
    \EE(W) \leq 2\binom{n}{k} (1-p)^{\binom{k}{2}}\sum_{r=1}^{n-k}\left[ 2\left(\frac{n}{k}\right) (1-p)^{k/2 - 1/2}(1-p/2+p^2)^{k+r} \binom{k+r}{2} \right]^{r}. \label{common_ratio}
\end{align}
Once again, one can show by the ratio test (with a little additional care) that the expression inside the square brackets will increase if $r$ is set to 0, which then becomes $n^{1 - 2(1 +\gamma) - \gamma +o(1)} = n^{ -1 - 3\gamma +o(1)} = o(1)$. Therefore,
\begin{align*}
    \EE(W) \le 2
    n^{1/3 + \gamma}\cdot n^{ -1 - 3\gamma +o(1)}
    = n^{ -2/3 - 2 \gamma + o(1)} = o(1).
\end{align*}
Thus, (ii) holds whp.

To verify (iii), we again define $W'$ to be the number of special subgraphs of order $k+1$. Adapting (\ref{common_ratio}) we have
\[ \EE[W'] \leq 2\binom{n}{k+1} (1-p)^{\binom{k+1}{2}}\sum_{r=0}^{n-k-1}\left[ 2\left(\frac{n}{k+1}\right) (1-p)^{k/2}(1-p/2+2p^2)^{k+1+r} \binom{k+1+r}{2} \right]^{r}. \]
We work with $ \binom{n}{k+1} (1-p)^{\binom{k+1}{2}}$ as before:
\begin{align*}
    \binom{n}{k+1} (1-p)^{\binom{k+1}{2}} &= \binom{n}{k} (1-p)^{\binom{k}{2}} \left(\frac{(1-p)^{k} (n-k)}{k+1}\right) 
    \\& \leq
    n^{\gamma+1/3} \cdot n^{ - 1 -\gamma+o(1)}
    \\&=
    n^{-2/3 +o(1)} \\&=
    o(1).
\end{align*}
As this term is $o(1)$, the $r=0$ term is $o(1)$, which implies that the sum of the remaining terms also vanishes (following the calculations from the previous case). This verifies that (iii) holds whp.
\end{proof}

Note that part i of Theorem~\ref{main} follows from Propositions \ref{prop_constant_p} and \ref{prop_nonconstant_p}: assume that (i), (ii), and (iii) from Propositions \ref{prop_constant_p} and \ref{prop_nonconstant_p} all hold. Then by (i) and (iii), $\a(G) = k-1$ or $k$. If $\a(G) = k-1$, then there exists no independent set of size $k$. By (ii), no other special subgraph of order $k$ exists either, hence by Lemma \ref{special}, $bp(G) \geq n - k + 1 = n - \a(G)$ as desired. If instead $\a(G) = k$, then by (iii) and Lemma \ref{special}, $bp(G) \geq n - k = n - \a(G)$, as desired.

The astute reader will notice that Proposition~\ref{prop_nonconstant_p} implies two-point concentration of the independence number 
of $ G_{n,p}$ for $ n^{-1/3 + \epsilon} <p < n^{ -\epsilon}$. Despite the fact that the proof of this fact follows the standard second moment
method approach to the independence number, we were not able to find this fact explicitly stated in the literature. So we state 
this as a Corollary of Proposition~\ref{prop_nonconstant_p}:
\begin{corollary}
If $ p = n^{\gamma}$ where $ -1/3 < \gamma < 0$, then there is a sequence $ k = k(n, \gamma)$ such that
$ \alpha( G_{n,p}) \in \{k-1,k\}$ whp. 
\end{corollary}

\section{Lower bound on $ bp( G_{n,p})$ for $p_0 \le p < 1/2$}

We prove part~ii of Theorem~\ref{main} using a first moment argument along the lines of the proof of part~i of Theorem~\ref{main}, but here the argument is more intricate and uses the BKR inequality. (See Section~\ref{bkr} for discussion of this inequality.)

Recall that we say that a biclique in $ G=G_{n,p}$ defined by vertex parts $B,B'$ is {\em based} on 
$B$ if $|B| \leq |B'|$ 
(so balanced bicliques are based on two sets).
We bound the probability that a particular set $A$ of $n'=k+r$ vertices in $ G_{n,p}$ forms a special subgraph of order $k$ by summing 
over the sets on which the bicliques that form the special subgraph can be based.
For  $A_1, \dots, A_r \subset A$ such that $|A_i| \ge 2$ for all $i$ and $ r < n' $, 
we define the event $ \cE= \cE_{A:A_1, \dots, A_r}$ to be the event
that the induced graph $G[A]$ is a special subgraph of order $k$ with the additional 
property that each biclique defining this special subgraph is based on at least one set from $A_1, \dots, A_r$. We begin by bounding the probability of such an event when we have the added 
condition that there is a uniform upper bound on the cardinality of each set $ A_i$. 
\begin{lem} \label{lenientprime}
Let $m \geq 2$ be a constant,
$A$ be a set of $n'$ vertices in $ G=G_{n,p}$, and $ A_1, \dots, A_r \subset A$ such that $ 2 \le |A_i| \le m$ for $i=1,\dots, r$, where $r<n'$.
If $n'$ is sufficiently large, then
\[ \PP( \cE_{A:A_1, \dots, A_r} ) \le (1-2^{-12}p)^{\frac{\binom{n'}{2}}{20 m^2}}. \]
\end{lem}
\begin{proof}
For $ j =2, \dots, m$, let $ r_j = |\{ i: |A_i| = j \}|$. So $r_j$ is the number of `bases' of size $j$ 
in the specified collection. 

Next, for any $A_{i}$, and any $v \not\in A_{i}$, define $$E'_{i}(v) = \left\{uv \in \binom{A}{2} \ \Bigg| \ u \in A_{i}\right\}.$$ Note that $|E'_{i}(v)| = |A_i|$, and that for every possible value of $j$, there are $(n'-j)r_j$ sets $E'_{i}(v)$ with size $j$. Note that the event $ \cE = \cE_{A:A_1, \dots, A_r}$ is contained in the event that the edge set of $G[A]$ can be partitioned into sets of the form $E'_i(v)$; call the latter event $\mathcal{E}'$.  Indeed, if $G[A]$ is decomposed into the bicliques defined by the sets
$ A_i, B_i$ for $i =1, \dots ,r$, then the edge set of $G[A]$ is the {\em disjoint} union of the sets in the collection $ \{ E'_{i}(v): i \in \{1, \dots, r\} \text{ and } v \in B_i \}$.

In order to bound $\PP[\mathcal{E}']$, 
define a hypergraph $\cL$ as follows: let $V(\cL) = \binom{A}{2}$ (so each vertex in $\cL$ corresponds to an unordered pair of vertices from $A$), and let the hyperedges in $\cL$ be the sets $E'_{i}(v)$. 
Note that we can view $ G[A]$ as being given by choosing
each vertex $v \in V(\cL)$ independently with probability $p$, and in order to prove the lemma it suffices to
bound the probability that there exist mutually disjoint edges $e_1,e_2,\dots,e_{\ell} \in \cL $ such that the set of chosen vertices (which
corresponds to the edge set of $G[A]$) is exactly $\bigcup_{i=1}^{\ell} e_i$, or in other words, there exists a subset of $E(\mathcal{L})$ that is a perfect matching of the random subset of $V(\mathcal{L})$.

For all $v \in V(\cL)$, $ j \in \{2, \dots, m\},$ let $N_j(v)$ be the number of $j$-edges of $ \cL$ that contain $v$. Then define two functions $f,g : V(\cL) \rightarrow \RR$:
\begin{align*}
    f(v) = \sum_{j=2}^{m} (j-1) N_j(v), \qquad
    g(v) = \sum_{j=2}^m  p^{j-1} N_j(v).
\end{align*}
Note that $f(v)$ gives an upper bound on the size of the neighborhood of $v$ in $ \cL$ if we define this neighborhood
to be the set of vertices $u$ such that $u$ and $v$ are contained in a common edge of $\cL$. It is straightforward to deduce that
\begin{align*}
    \sum_{v \in V(L)} f(v) < \sum_{j=2}^{m}j^2 r_j n' \leq m^2 n' \sum_{j=2}^m r_j \leq m^2n'(n'-1)
\end{align*}
and
\begin{align*}
    \sum_{v \in V(L)} g(v) < \sum_{j=2}^{m} (j p^{j-1}) r_j  n'< n'\sum_{j=2}^{m} r_j \leq n'(n'-1) \qquad \text{(since $p < 1/2$).}
\end{align*}
It follows that at least $ \frac{2}{3} \binom{n'}{2}$ vertices in $\cL$ have an $f$-value of at most $6m^2$, and at least $ \frac{2}{3} \binom{n'}{2}$ vertices in $\cL$ have a $g$-value of at most $6$. Therefore at least $ \frac{1}{3} \binom{n'}{2}$ vertices have both properties. Let $S$ be the set of such vertices.

Next, we construct a set $T \subseteq S$ with the property that there is no edge of $\cL$ that contains two elements of $T$. We construct
this set by iteratively choosing elements of $S$ to be elements of $T$ one at a time, deleting all vertices that are in a common edge with
the selected vertex at each stage.  We have
$|T| \geq |S|/(6m^2 + 1) \geq \binom{n'}{2}/(20m^2)$.

We now fix such $T$ and consider the probability space. We reveal the choice of elements of $ \cL$ (which corresponds to the choice of edges in $ G[A]$) in two stages. First, we reveal the elements of $T$ that are chosen; let this set be $I$. We then reveal the vertices chosen among $V(\cL) \backslash T$; let this set be $I'$. Now, given the set $I$, the event $\cE$ is contained in the event that for all $v \in I$, there exists some $e_v \in \cL$ such that $e_v \ni v$ and $e \subset I \cup I'$.  Furthermore, and crucially, these $|I|$ edges of $ \cL$ are pairwise disjoint. 
Indeed, each element $v \in I$ corresponds to an edge of $G[A]$, and once chosen, this edge must be included in one of the 
bicliques. From the perspective of this fixed $ v \in V(\cL)$, any $ e \ni v$ such that $ e\setminus v \subset I'$ gives a potential
choice for this biclique. However, as the edge sets of the bicliques are pairwise disjoint, 
the collection of choices $( e_v : v \in I )$ must be made so that these edges of $\cL$ are disjoint. So, let $C$ be the
event that there is a collection $ (e_v: v \in I)$ of edges of $ \cL$ with the following properties:
\begin{itemize}
    \item[(a)] $ v \in e_v$ for all $ v \in I$,
    \item[(b)] $ e_v \setminus v \subseteq I'$ for all $ v \in I$, and
    \item[(c)] if $ u \neq v$, then $ e_u \cap e_v = \emptyset$.
\end{itemize}
We have
\begin{align*}
    \PP[\cE] \leq  \PP[ \cE'] \leq \sum_{I \subseteq T} p^{|I|}(1-p)^{|T|-|I|}\PP[C| I].
\end{align*}
It remains to bound $\PP[C | I]$.
For each $v \in I$, and each $e \ni v$, let $D_{v,e}$ be the event that $e \setminus v \subset I'$, so $\PP[D_{v,e}] = p^{|e|-1}$.
Moreover, for a fixed $v$ the events $D_{v,e}$ are increasing, and therefore by FKG (e.g. the second part of Theorem~6.3.2 in \cite{probmethod}) we have
\begin{align*}
    \PP\left[\bigcap_{e \ni v} \overline{D}_{v,e}\right] \geq \prod_{e \ni v}\PP[\overline{D}_{v,e}] = \prod_{j=2}^{m}(1-p^{j-1})^{N_j(v)}.
\end{align*}
Since $p < 1/2$, and by using the inequality $\ln(1-x) > -2 \ln(2)x$ for $x < 1/2$ (derivable from elementary calculus), it follows that
\begin{multline*}
    \ln\left(\PP\left[\bigcap_{e \ni v} \overline{D}_{v,e}\right]\right) \geq \ln\left( \prod_{j=2}^{m}(1-p^{j-1})^{N_j(v)}\right) = \sum_{j=2}^{m} N_j(v) \ln(1-p^{j-1}) \\
    > \sum_{j=2}^{m} N_j(v) (-2 \ln (2) p^{j-1}) = - 2 \ln (2) g(v) \geq - 12 \ln(2).
\end{multline*}
Hence, $\PP\left[\bigcap_{e \ni v} \overline{D}_{v,e}\right] > 2^{-12}$.

Now we are ready for our application of the BKR inequality (see Section~\ref{bkr} for definitions and discussion). Define $C_v = \bigcup_{e \ni v} D_{v,e}$.  Note that a cylinder $ [C_v]_J$ is nonempty if and only if 
$J$ contains a set $e \setminus v$ for some $ e \in \cL$ that contains $v$. Thus, conditioning on $I$, we have
\[ C = \square_{v \in I} C_v\]
and, applying the BKR inequality (to be precise, Theorem~\ref{BKRmany}), we have
\begin{align*}
    \PP[C | I] \leq \prod_{v \in I} \PP[C_v] < (1 - 2^{-12})^{|I|}.
\end{align*}
We complete the proof of the lemma with an application of the binomial theorem:
\begin{align*}
    \PP( \cE_{A; A_1, A_2, \dots, A_r}) &\leq \sum_{I \subseteq T} p^{|I|}(1-p)^{|T|-|I|}\PP[C| I] \\&\leq
    \sum_{I \subseteq T} (p(1-2^{-12}))^{|I|}(1-p)^{|T|-|I|} \\&=
    (1 - 2^{-12}p)^{|T|} \\&\leq
    (1 - 2^{-12}p)^{\frac{\binom{n'}{2}}{20 m^2}}.
\end{align*}
\end{proof}

Next we account for the contribution to $ {\mathbb P} ( {\mathcal E}_{A: A_1, \dots, A_r})$ 
coming from bicliques based on sets $ A_i$ with $|A_i| >m$. This is achieved in the following
two lemmas.

\begin{lem} \label{lem: include big m}
Let $A$ be a set of $n'$ vertices in $G = G_{n,p}$, $s > s'$ be positive integers, and $A_1, \dots, A_s \subseteq A$ such that $|A_i| \geq 2$ for all $i$. Then
    $$\PP( \cE_{A; A_1, \dots, A_s }) \le \PP( \cE_{A; A_1, \dots, A_{s'} }) \sum_{ B_{s'+1} ,\dots, B_s} 
\prod_{j=s'+1}^s  \left( \frac{ p}{ 1-p} \right)^{ |A_j| \cdot |B_j|}$$
where the sum is over all subsets $ B_{s'+1}, \dots, B_s$ of $A$ such that 
$ |B_j| \ge |A_j| $ and $A_j \cap B_j = \emptyset$ for $ j =s'+1, \dots, s$.
\end{lem}

\begin{proof}
     First, fix partite set pairs $(A_i,B_i)$ for $i \in [s]$ such that $A_i \cap B_i = \emptyset$. The probability of $G[A]$ being the disjoint union of the $s$ bicliques with these partite sets is at most $${(1-p)^{\binom{n'}{2}} \prod_{j=1}^{s} \left(\frac{p}{1-p}\right)^{|A_j||B_j|}}$$ (this holds with equality if the bicliques are disjoint; the probability is 0 otherwise). Next, note that, if one sums the previous expression over all possible $B_1, \dots, B_{s'}$ such that $|B_i| \geq |A_i|$ for $i \in [s']$, one arrives at the expression $$\PP( \cE_{A; A_1, \dots, A_{s'} }) \prod_{j=s'+1}^s  \left( \frac{ p}{ 1-p} \right)^{ |A_j| \cdot |B_j|}$$ by definition of $\mathcal{E}$. Finally, we allow the rest of the $B_i$ (with $i > s'$) to vary; by summing over all such $B_i$ with $|B_i| \geq |A_i|$, one arrives at the desired inequality.
\end{proof}

\color{black}

\begin{lem} \label{modified 2}
Let $m$ be a constant,
and $A$ be a set of $n'$ vertices in $ G=G_{n,p}$.
If $n'$ is sufficiently large (as a function of $m$ and $p$), then
the probability $P_{A:r}$ that $G[A]$ can be partitioned into at most 
$r$ nonempty non-star bicliques is at most
\begin{align}
    4(1- 2^{-12}p)^{\frac{\binom{n'}{2}}{20m^2}} \left(2 \binom{n'}{m+1} \left(1+\left(\frac{p}{1-p}\right)^{m+1}\right)^{n'}\right)^{r} \label{combine to r}.
\end{align}
\end{lem}
\begin{proof}
Let $ P_{A,s}$ be the probability that $ G[A]$ can be partitioned into exactly $s$ nonempty non-star bicliques.
We bound $ P_{A,s}$ using the union bound, summing over all possible choices of base sets $ A_1, \dots, A_s \subset A$.  Suppose such a collection of
sets is fixed, and further suppose that  $ |A_i| \le m$ for $ i =1, \dots, s'$ and $|A_i| > m$ for $ i = s'+1, \dots, s$. By Lemma \ref{lem: include big m}, 
we have
\[ \PP( \cE_{A; A_1, \dots, A_s }) \le \PP( \cE_{A; A_1, \dots, A_{s'} }) \sum_{ B_{s'+1} ,\dots, B_s}  
\prod_{j=s'+1}^s  \left( \frac{ p}{ 1-p} \right)^{ |A_j| \cdot |B_j|},   \]
where the sum is over all subsets $ B_{s'+1}, \dots, B_s$ of $A$ such that 
$ |B_j| \ge |A_j| $ and $A_j \cap B_j = \emptyset$ for $ j =s'+1, \dots, s$. With this inequality in hand, and applying Lemmas~\ref{combo}~and~\ref{lenientprime} to get (\ref{equ: tight for space}) below, we can write

\begin{align}
    P_{A,s} & = \sum_{s' = 0}^s \sum_{A_1, \dots, A_s' \in \binom{A}{ \le m}} \sum_{A_{s'+1}, \dots, A_s \in \binom{A}{ >m}  }\PP( \cE_{A; A_1, \dots, A_s}) \nonumber \\
    & \leq  \sum_{s' = 0}^s \sum_{A_1, \dots, A_s' \in \binom{A}{ \le m}}  \PP( \cE_{A; A_1, \dots, A_{s'}}) 
    \left( \sum_{\l \geq \mu \geq m+1, \l + \mu \leq n'} \binom{n'}{\l + \mu} \binom{\l + \mu}{\l} \left(\frac{p}{1-p}\right)^{\l \mu} \right)^{ s - s'} \nonumber \\
    & \leq  \sum_{s'=0}^s \binom{n'}{m+1}^{s'} (1- 2^{-12}p)^{\frac{\binom{n'}{2}}{20m^2}}  
     \left( 2 \binom{n'}{m+1} \left( 1 + \left( \frac{p}{1-p}\right)^{m+1} \right)^{n'}\right)^{ s - s'} \label{equ: tight for space} \\
    & = (1- 2^{-12}p)^{\frac{\binom{n'}{2}}{20m^2}}  
     \binom{n'}{m+1}^s \sum_{s'=0}^s  \left( 2 \left( 1 + \left( \frac{p}{1-p}\right)^{m+1} \right)^{n'}\right)^{ s - s' } \nonumber \\
     & \leq 2  (1- 2^{-12}p)^{\frac{\binom{n'}{2}}{20m^2}}  
     \left( 2 \binom{n'}{m+1}  \left( 1 + \left( \frac{p}{1-p}\right)^{m+1} \right)^{n'}\right)^{ s}, \nonumber
\end{align}
where, in the last step, we use the fact that the sum is a geometric series with common ratio less than 1/2. Summing over all the possible values of $s$ gives the bound in the lemma, again observing that we have a geometric series.
\end{proof}

We are now ready to complete the proof of part ii of Theorem \ref{main}. Recall that $p$ is a constant such that $ p_0 \le p < 1/2$. Let $W$ be the random variable defined as the number of nonempty (edge-wise) special subgraphs of order $k$ in $G_{n,p}$, where $k$ will be chosen later 
to be $O(\log(n))$. By Lemma \ref{special} it is enough to show that $\EE(W) = o(1)$. 

We apply Lemma~\ref{modified 2} with $n' = k + r$, where $m = m_p$ is defined to be an integer such that $m_p > 3$ and
\begin{align*}
    (1-2^{-12}p)^{\frac{-1}{80 m_p^2}} > 1+\left(\frac{p}{1-p}\right)^{m_p+1};
\end{align*}
such an $m$ exists since the left side is $1+\Theta(m^{-2})$ and the right side is $1+e^{-\Theta(m)}$. For sufficiently large $n$ and $k$, and recalling that both $p$ and $m_p$ are constants, we have
\begin{align*}
    \EE(W) &\leq  \sum_{r = 0}^{\infty}  \sum_{A \in \binom{[n]}{k+r}} P_{A:r} \\
    & \leq \sum_{r=0}^{\infty}
    \binom{n}{k+r} 4(1- 2^{-12}p)^{\frac{\binom{k+r}{2}}{20 m_p^2}} \left(2 \binom{k+r}{m_p+1} \left(1+\left(\frac{p}{1-p}\right)^{m_p+1}\right)^{k+r}\right)^{r} \\
    &\leq 
    4 \sum_{r=0}^{\infty} n^{k+r}(1- 2^{-12}p)^{\frac{(k+r)^2}{40m_p^2}} \left(1+\left(\frac{p}{1-p}\right)^{m_p+1}\right)^{2(k+r)r} \\&\leq
    4 \sum_{r=0}^{\infty} n^{k + r} (1- 2^{-12}p)^{\frac{(k+r)k}{40m_p^2}} \\&=
    4\sum_{r=0}^{\infty}  \left(n (1- 2^{-12}p)^{\frac{k}{40m_p^2}}\right)^{k+r}.
\end{align*}
Now choose $k$ such that 
\begin{align*}
    (1-2^{-12}p)^{\frac{k}{40 m_p^2}} < 1/n^2.
\end{align*}
This can be done with some $k = O(\ln(n))$. This will result in $\EE(W) = o(1)$.

\section{Upper bound on $bp(G_{n,p})$ for $ p_0 < p < 1/2 $}

Here we prove part iii of Theorem \ref{main}, using the second moment method and a variation of the proof of the main result in \cite{ABHmore}. Throughout this section we will omit all floor and ceiling signs and ignore divisibility issues whenever these are not crucial, so as to simplify the presentation. For ease of notation we set $b = 1/(1-p)$.

First, recall that by Lemma \ref{root}:
\begin{align*}
    (1 - p)^{-\frac{1}{2}}  < 1 + \left( \frac{p}{1 - p}\right)^2.
\end{align*}
Furthermore, one can check that, since $p < 0.5$, then $(1-p)^2 + p^2 < 1-p$.
Next, let $a_p$ and $c_p$ be  positive constants such that $a_p < 0.01$, $c_p < 0.01$, and such that the following three inequalities hold (think of setting $c_p \ll a_p \ll 1$):
\begin{align}
    \left[(1 - p)^{-\frac{1}{2}} \right]^{ \frac{a_p + 4 c_p}{a_p (1-a_p)(1+2c_p)}} &< 1 + \left( \frac{p}{1 - p}\right)^2  \label{possible}\\
    ((1-p)^2 + p^2)^{\frac{1-10a_p}{1-a_p/2}} &< 1-p \label{Possible}\\
    \frac{(1-4a_p)(1+c_p)}{1-a_p/2} &< 1. \label{easily}
\end{align}

Now set 
\begin{align*}
    k  = 2 (1+c_p)\log_b(n)/(1-a_p/2) 
\end{align*}
(note that this is slightly different than what $k$ was before). Define a family $\mathcal{F}_k'$ of graphs on $k$ vertices as follows: a graph $F$ with $|F| = k$ is in $\mathcal{F}_k'$ if and only if the following two conditions hold: 
\begin{enumerate}
    \item $F$ is a bipartite graph with partite sets $A$ and $B$ such that $|A| = a_p k$ and $|B| = (1 - a_p)k$.
    \item Set $A$ is partitioned into sets $A_1,A_2,\dots,A_r$, where $r = a_p k/2$ and $|A_i| = 2$ for all $i \in [r]$. Moreover, for every $A_i$ and $\b \in B$, either $\b$ is adjacent to both vertices of $A_i$ or neither of them.
\end{enumerate}
Next, let $\mathcal{F}_k$ be a subfamily of $\mathcal{F}_k'$ such that the two \textit{additional} conditions hold for all $F \in \mathcal{F}_k$: 
\begin{enumerate}
    \item[3.] The degree of each $a \in A$ is at least $0.15k$. \label{cond 3}
    \item[4.] For each $i,j$ corresponding to different sets $A_i,A_j$, the symmetric difference of the neighborhoods of $A_i$ and $A_j$ is at least $0.25k$. \label{cond 4}
\end{enumerate}
Note that if there exists some induced subgraph $F \subseteq G_{n,p}$ such that $F \in \mathcal{F}_k'$, then we get the desired result: $E(F)$ is a disjoint union of (at most) $r = a_p k/2$ bicliques (where each of these have an $A_i$ for one partite set), and there are $k$ vertices in total. Hence, applying a well-known bound on the independence number of the dense random graph (see Section~7.2 of \cite{alanmichal}), we have
\begin{multline*}
    bp(G_{n,p}) \leq n - k + r = n - (1 - a_p/2)k \\ = n - 2(1+c_p)\log_b(n) \leq n - (1+c_p)\a(G_{n,p}) \qquad  \text{whp.}
\end{multline*}

Next, we introduce notation to help us define a suitable random variable. For any set $K \subset V(G)$ with $|K| = k$, let $M_K$ be the event that the subgraph induced on $K$ is in $\mathcal{F}_k$. We then let $X_K$ be the indicator random variable of $M_K$, and finally, let $X = \sum_{K} X_K$. Hence, $X$ is the number of induced subgraphs of $G_{n,p}$ that are in $\mathcal{F}_k$. We want to show that $X \geq 1$ whp using the second moment method; specifically, we show that $\EE(X) = \omega(1)$ and that $\text{Var}(X) = o(\EE(X)^2)$, which together with Chebyshev's inequality imply that $X \geq 1$ whp. We emphasize that, though we only need to show the existence of a subgraph of $G_{n,p}$ in $\mathcal{F}_k'$ whp, the two additional conditions of $\mathcal{F}_k$ make the second moment calculations manageable.

\subsection{First Moment Calculation}
To start, we introduce a (slightly) different random variable $X'$ that equals the number of induced subgraphs of $G_{n,p}$ which are in $\mathcal{F}_k'$. We now seek to find a lower bound for $\EE(X')$. To start, first calculate the number of ways to choose a set of $k$ vertices, and among them choose partite sets $A$ and $B$, and furthermore split the set $A$ into $r$ sets of size 2. This number is
\begin{align}
    \binom{n}{k} \binom{k}{2} \binom{k-2}{2} \dots \binom{k-2r+2}{2} \frac{1}{r!} \nonumber &= \binom{n}{k} \frac{k!}{2^{r} (k-2r)!r!} \\&=
    \binom{n}{k} \frac{k!}{2^{a_p k/2}((1-a_p)k)!(a_p k/2)!} \nonumber
    \\&=
    n^k k^{(-1 + a_p/2 \pm o(1))k}.  \nonumber
\end{align}
Next, we bound below the probability of there being an induced graph in $\mathcal{F}_k'$ given $A_1,\dots,A_r,B$. For each pair consisting of a vertex $\b \in B$ and a set $A_i$, the probability of there being either two edges or two non-edges between $\b$ and $A_i$ is $(1-p)^{2} + p^{2}$. There are $(1 - a_p)kr$ such pairs, so the probability of picking edges accordingly from $A$ to $B$ is $((1-p)^{2} + p^{2})^{(1 - a_p)kr}$. Hence, we have the following lower bound on the total probability of a graph $F \in \mathcal{F}_k'$ being produced, 
given $A_1,A_2,\dots,A_r,B$:
\begin{align*}
    & \phantom{=} ((1-p)^{2} + p^{2})^{(1-a_p)kr} (1-p)^{\binom{(1 - a_p)k}{2} + \binom{a_p k}{2}}
    \\ &>
    ((1-p)^{2} + p^{2})^{(1-a_p)kr}  (1-p)^{(a_p^2 - a_p + 1/2)k^2} \\&=
    (1 + p^{2}/(1-p)^2)^{(1-a_p)kr}  (1-p)^{(a_p^2 - a_p + 1/2)k^2 + 2(1 - a_p)kr}
    \\&=
    (1 + p^{2}/(1-p)^2)^{((a_p-a_p^2)/2) k^2}  (1-p)^{k^2/2} \\&=
    \left((1 + p^{2}/(1-p)^2)^{((a_p-a_p^2)/2)}(1-p)^{1/2} \right)^{k^2}.
\end{align*}
Putting these together, we have
\begin{align}
    \EE(X') > n^k k^{(-1 + a_p/2 \pm o(1))k}\left((1 + p^{2}/(1-p)^2)^{((a_p-a_p^2)/2)}(1-p)^{1/2} \right)^{k^2}. \label{equ: exp X' bound}
\end{align}
The following Lemma will allow us to bound $\EE(X)$ by showing that a typical graph in $\mathcal{F}_k'$ is also in $\mathcal{F}_k$:
\begin{lem} \label{lem: X' to X}
    Take a vertex set $K$ with size $k$ and fixed partitions $ K = A \cup B$ and $ A = A_1 \cup \dots \cup A_r$ with $|A| = a_p k$, $ |A_i|=2$ for all $i$, and  $|B| = (1-a_p)k$. If we take $G_{k,p}$ on set $K$ and condition on $G_{k,p} \in \mathcal{F}_k'$ with respect to our given partitions, then whp $G_{k,p} \in \mathcal{F}_k$.
\end{lem}
\begin{proof} To start, note that, for any pair $(A_i,\b)$ with $ \b \in B$, the conditional probability that $\b$ is adjacent to the vertices in $A_i$ equals $p^2/( (1-p)^2 + p^2) $, and these events are independent across all $\b \in B$ and all $A_i$. One can also check that, for all $p \in (p_0,0.5)$, that  $ p^2/( p^2 + (1-p)^2) \in (0.17,0.5)$. Therefore, the degree of any $a \in A$ is a binomial random variable with expectation $(1 - a_p)k p^2/( p^2 + (1-p)^2) $, which is at least $0.16k$ (since $a_p < 0.01$). Therefore, one can use a standard Chernoff bound (e.g. equation (2.9) in \cite{janson2011random}), to deduce that the probability that a specific vertex in $A$ has degree less than $0.15k$ is $e^{-\Omega(k)}$, so condition 3 is satisfied whp by using a union bound. Similar logic applies to condition 4: the symmetric difference of the neighborhoods of $A_i$ and $A_j$ for $i \not= j$ is a binomial random variable with expectation 
$ (1-a_p)k \cdot 2p^2(1-p)^2/((1-p)^2 + p^2)^2$. One can check that for $p \in (p_0,0.5)$ we have $ 2p^2(1-p)^2/ ((1-p)^2 + p^2)^2 \in (0.28,0.5)$. Hence, the expected symmetric difference will be at least $0.27k$, and by using Chernoff and union bounds once again, condition 4 is satisfied whp.
\end{proof}

By (\ref{equ: exp X' bound}) and Lemma \ref{lem: X' to X}, \color{black} we have
\begin{align*}
    \EE(X) &> n^k k^{(-1 + a_p/2 \pm o(1))k}\left((1 + p^{2}/(1-p)^2)^{((a_p-a_p^2)/2)}(1-p)^{1/2} \right)^{k^2} \\&=
    \left(n k^{-1 + a_p/2 \pm o(1)} \bigg((1 + p^{2}/(1-p)^2)^{((a_p-a_p^2)/2)}(1-p)^{1/2} \bigg)^k \right)^k.
\end{align*}
By (\ref{possible}), we have that
\begin{align*}
    (1 + p^{2}/(1-p)^2)^{(a_p-a_p^2)/2} > (1-p)^{-\frac{a_p+4c_p}{4(1+2c_p)}}.
\end{align*}
Therefore
\begin{align*}
    \EE(X) > \left(n k^{-1 + a_p/2 \pm o(1)}\bigg((1-p)^{\frac{2 - a_p}{4+8c_p}} \bigg)^k \right)^k.
\end{align*}
We want to show that $\EE(X) = \omega(1)$, hence it is enough to show that 
\begin{align}
    \frac{n}{\ln(n)} (1-p)^{\frac{(2 - a_p)k}{4+8c_p}} = \omega(1). \label{plug in k}
\end{align}
Now we plug in our value of $k$ into the left side of (\ref{plug in k}) to get
\begin{align*}
    & \phantom{=} \frac{n}{\ln(n)} (1-p)^{\frac{(2 - a_p)(2 (1+c_p))\log_b(n)}{(4+8c_p)(1-a_p/2)}} =
    \frac{n^{c_p/(1+2c_p)}}{\ln(n)} = \omega(1).
\end{align*}
Therefore, $\EE(X) = \omega(1)$ as desired.

\subsection{Second Moment Calculation}
We now proceed to estimate Var$(X)$. Similarly to the proof of Proposition~\ref{prop_constant_p}, for  $K,K' \subset V(G),|K|=|K'|=k$, let $K \sim K'$ denote that $| K \cap K'|\geq  2$ (and $K \not= K'$). The variance of $X$ satisfies
\begin{align*}
    \text{Var}(X) = \sum_K \text{Var}(X_K)+ \sum_{K \sim K'} \text{Cov}(X_K,X_{K'}) \leq \EE(X)+ \sum_{K \sim K'} \EE(X_K X_{K'}).
\end{align*}

For each $i$, $2 \leq i \leq k - 1$, let $h_i$ denote the contribution of the pairs with intersection size $i$ to the
above sum, that is:
\begin{align*}
    h_i = \sum_{|K \cap K'| = i} \EE(X_K X_{K'}).
\end{align*}
(Note that this is slightly different than the notation in Proposition \ref{prop_constant_p}.) 
Our objective is to show that $\sum_{i=2}^{k-1} h_i = o(\EE(X)^2)$; this, combined with $\EE(X) = \omega(1)$, gives us the desired result.
We consider two possible ranges for the parameter $i$ as follows:

\noindent
{\bf Case 1:} $i \leq (1 - 10 a_p)k$

\noindent
Here we show that $(\sum_{i=2}^{(1-10 a_p)k}$ $h_i)/\EE(X)^2 = o(1)$. Recall that $M_K$ is the event that the induced subgraph $F$ on vertex set $K$ is in $\mathcal{F}_k$. We have $\binom{n}{k} \binom{k}{i} \binom{n-k}{k-i}$ ways to choose sets $K$ and $K'$ such that $|K \cap K'| = i$, so
\begin{align*}
    h_i = \binom{n}{k} \binom{k}{i} \binom{n-k}{k-i} \PP[M_K \cap M_{K'}].
\end{align*}

\begin{lem}
For all $i \leq (1-10 a_p)k$ with $K,K' \in \binom{[n]}{k}$ and $|K\cap K'|=i$, we have
\begin{align*}
    \PP[M_K \cap M_{K'}] \leq (1-p)^{-(0.5 - 2a_p)ki} (\PP[M_K])^2.
\end{align*}
\end{lem}
\begin{proof}
Let $F \in \mathcal{F}_k$, and let $F'' \subseteq F$ be a subgraph with $|F''| = i$. ($F''$ represents the subgraph induced by $K \cap K'$.) Then since $F''$ is bipartite, with one partite set with at most $a_p k$ vertices, then $F''$ must have less than $a_p k i$ edges. Also, it can be quickly checked that, for all $p \in (p_0,0.5)$, that $(1-p)^4 < p < 1-p$. It follows that
\begin{align*}
    \PP[M_K \cap M_{K'}] &< (1-p)^{-\binom{i}{2}+a_p k i} p^{-a_p k i} (\PP[M_K])^2 \\&<
    (1-p)^{-i^2/2 + a_p k i} (1-p)^{-4 a_p k i} (\PP[M_K])^2 \\&=
    \big((1 - p)^{-i/2 - 3a_p k}\big)^i (\PP[M_K])^2 \\&\leq
    (1-p)^{-(0.5 - 2a_p)ki}(\PP[M_K])^2.
\end{align*}
\end{proof}

Since $\EE(X) = \binom{n}{k}\PP[M_K]$, we have
\begin{align}
    \frac{h_i}{\EE(X)^2} \leq \frac{\binom{k}{i} \binom{n-k}{k-i}}{\binom{n}{k}}(1-p)^{-(0.5 - 2a_p)ki} &\leq k^i \left(\frac{k}{n} \right)^i (1-p)^{-(0.5 - 2a_p)ki} \nonumber \\&= \left(\frac{k^2 b^{(0.5 - 2a_p)k}}{n} \right)^i \nonumber \\&=
    \left(k^2 n^{\frac{(1 - 4a_p)(1+c_p)}{1-a_p/2}-1}\right)^i. \label{2.2.1}
\end{align}
By (\ref{easily}), we have $\frac{(1-4a_p)(1+c_p)}{1-a_p/2} < 1$, and thus
\begin{align*}
    k^2 n^{\frac{(1-4a_p)(1+c_p)}{1-a_p/2}-1} = o(1).
\end{align*}
Summing the expression in (\ref{2.2.1}) over $i \leq (1-10a_p)k$ gives us the desired result.

\noindent
{\bf Case 2:} $i > (1-10a_p)k$

\noindent
We now show that $(\sum_{i=(1-10a_p)k}^{k-1}$ $h_i)/\EE(X)^2 = o(1)$ by a different method.
First, we use the secondary conditions to impose restrictions on the intersections of two $\mathcal{F}_k$ graphs in this case.
\begin{lem}
If $F,F' \in \mathcal{F}_k$ such that $|V(F \cap F')| = |K \cap K'| \geq (1-10a_p)k$, then the following hold:
\begin{enumerate}
    \item A vertex $v \in K \cap K'$ is in partite set $A$ of $K$ if and only if $v$ is in partite set $A'$ of $K'$.
    \item If $u,v \in A \cap A'$, then $u,v$ are paired together in some set $A_{j}$ of $A$ if and only if they are paired together in some set $A_{j}'$ of $A'$.
\end{enumerate}
\label{intersection}
\end{lem}
\begin{proof} Since we set $a_p < 0.01$, we have $(1 - 10a_p)k > 0.9k$. Claim 1 of the lemma follows from condition 3 in the definition of $ {\mathcal F}_k$: if $v \in K \cap K'$ and $v \in A$, then the degree of $v$ in graph $F \cap F'$ is at least $0.15k - 0.1k > a_p k$. Since the degree of any $\b \in B'$ is at most $a_p k$, then $v \in A'$ also. The converse follows by symmetry.

Claim 2 of the lemma follows from condition 4 in the definition of ${\mathcal F}_k$: if $u,v \in A \cap A'$ are {\it not} paired in the same set $A_{j}$ of $A$, then the neighborhoods of $u$ and $v$ within $F \cap F'$ still differ by at least $0.25k - 0.1k > 0$ vertices. Hence, they have different neighborhoods in $F'$, so they cannot be paired together in the same set $A_{j}'$ in $A'$. The converse follows by symmetry.
\end{proof}

Next, we find a bound on $\PP[M_K \cap M_{K'}]$; more precisely, a bound on $\PP[M_{K'} | M_{K}]$. We fix some $F \in \mathcal{F}_k$ on vertex set $K$, and bound  $\PP[M_{K'}]$. Lemma~\ref{intersection} tells us exactly how vertices from $K \cap K'$ should be placed into sets $A_1',\dots,A_r',$ and $B'$ (without loss of generality). Hence, our strategy is as follows: bound the number of ways to assign each vertex in $K'\backslash K$ to one of the $r+1$ sets $A_1',\dots,A_r',B'$, and then bound the probability of randomly picking edges and non-edges correctly for the $\binom{k}{2} - \binom{i}{2}$ remaining sites in $\binom{K'}{2}$ so that we have a graph in $\mathcal{F}_k$ with sets $A_1',\dots,A_r',B'$. For the first bound, we use the trivial expression  $(r+1)^{k-i}$. After such locations for the vertices are made, consider the probability that edges are picked or not picked in accordance with $M_{K'}$. Out of the $\binom{k}{2} - \binom{i}{2}$ sites, there may be certain pairs of sites that can either be both edges or both non-edges (with total probability $(1-p)^2 + p^2$). The rest of the sites are ``determined" (so we know whether each should contain or not contain an edge). For each of these sites, the probability that we choose (or not choose) an edge accordingly is at most $1-p$, which is bounded above 
by $((1-p)^2 + p^2)^{1/2}$. Combining this all together, the probability of the remaining edges being picked or not picked in accordance with $M_{K'}$, with the fixed $A_1', \dots, A_r',$ and $B'$, is at most $((1-p)^2 + p^2)^{\left(\binom{k}{2} - \binom{i}{2}\right)/2}$, which is itself bounded by $((1-p)^2 + p^2)^{i(k-i)/2}$.  Therefore,
\begin{align*}
    \PP[M_{K'} | M_K] &\leq (r+1)^{k-i} ((1-p)^{2} + p^{2})^{i(k-i)/2} \\&<
    \left((r+1) ((1-p)^{2} + p^{2})^{(1-10a_p)k/2}\right)^{k-i}.
\end{align*}

Finally, we bound $h_i / \EE(X)^2$:
\begin{align}
    \frac{h_i}{\EE(X)^2} < \frac{h_i}{\EE(X)}&=  \frac{\binom{n}{k} \binom{k}{i} \binom{n-k}{k-i} \PP[M_K \cap M_{K'}] }{\EE(X)} \nonumber \\&=
    \frac{\binom{n}{k} \binom{k}{i} \binom{n-k}{k-i} \PP[M_K]\PP[M_{K'} | M_{K}]}{\EE(X)} \nonumber \\&=
    \binom{k}{k-i} \binom{n-k}{k-i}\PP[M_{K'} | M_{K}] \nonumber \\& <
    \binom{k}{k-i} \binom{n-k}{k-i} \left((r+1) ((1-p)^{2} + p^{2})^{(1-10a_p)k/2}\right)^{k-i} \nonumber \\&\leq
    \left(kn(r+1)((1-p)^{2} + p^{2})^{((1-10a_p)/2)k} \right)^{k-i}. \label{almost there}
\end{align}
Now we plug in $k = \frac{2(1+c_p)\log_b(n)}{1-a_p/2}$ once again, and we use (\ref{Possible}) to bound (\ref{almost there}):
\begin{align}
    \frac{h_i}{\EE(X)^2} &\leq \left(kn(r+1)((1-p)^{2} + p^{2})^{((1-10a_p)/2)k} \right)^{k-i} \nonumber \\&\leq
    \left(kn(r+1)(1-p)^{\log_b(n)(1+c_p)}\right)^{k-i} \nonumber \\&=
    \left(k(r+1)n^{-c_p}\right)^{k-i}. \label{the end}
\end{align}
Clearly, $k(r+1)n^{- c_p} = o(1)$, so summing the expression in (\ref{the end}) over $i \in ((1-10a_p)k,k)$ gives us the desired result.

\section{Conclusion}

The main contributions of this work include establishing the critical probability $p_0$ for the property that $ bp( G_{n,p}) = n - \a ( G_{n,p})$ and 
establishing the order of magnitude of $n - bp(G_{n,p})$ for all constants $ p_0 \leq p < 1/2$. A number of interesting questions
remain open:
\begin{itemize}
\item It seems that the most interesting open question here is determining whether or not $ bp( G_{n,1/2})$ is bounded 
below by $ n - O( \log n)$. Alon conjectures that  $ bp( G_{n, 1/2}) > n - O( \log n)$ whp (this is Conjecture~4.1 in \cite{Alonbp}). The methods we introduce in Section~4 do not fully extend to $p=1/2$.  Very roughly speaking, we do not
have a way to suitably bound the probability of the appearance of bicliques in a special subgraph in which {\em both} of the vertex sets that define the
biclique are large. The lack of a way to handle this situation is reflected in the constant in part (ii) of Theorem~\ref{main}: as $p$ tends to $ 1/2$ from
below, we need to allow larger `bases' for these bicliques in our application of Lemma~\ref{lenientprime} (this is quantified by the constant $m_p$).
As a consequence, the constant in part (ii) of Theorem~\ref{main} goes to infinity as $ p$ tends to $1/2$ from below. Indeed, if we consider 
$ p = 1/2 - \epsilon$, then the constant in part (ii) of
Theorem~\ref{main} is on the order of $ (( 1/\epsilon) \log( 1/ \epsilon))^2$ as $ \epsilon \to 0$.

\item
The behavior of $ n - bp( G_{n,p} )$ for $p$ at and slightly below and above $p_0$ is not clear; all that is currently known 
is that $n - bp(G_{n,p}) = O(\ln(n))$ whp in this regime. 
It seems reasonable to speculate that $ bp( G_{n, p_0}) = n - \a( G_{n,p_0})$ whp and that there exists some smooth function $f$ such that $ \lim_{x \to 0^+} f(x) = 1$ and  $  n - bp( G_{n, p}) < f(\delta) \alpha( G_{n,p})$ whp for $ p = p_0 + \delta$ and $ \delta>0$ small.

\item
The behavior of $ n - bp( G_{n,p})$ is also not fully understood for $ p = n^{\gamma}$ where $ \gamma \in (-1,-1/3]$.  Alon showed that 
$  n - bp( G_{n,p}) = O( \a( G_{n,p}))$ in this regime and gave a precise characterization of 
$n - bp( G_{n,p})$ when $ p = n^{\gamma}$ where $ -1 < \gamma < -7/8$.  See Theorem~1.2 and Proposition~1.3 in \cite{Alonbp}.

\end{itemize}

For $ p > 1/2$ the question of bounding $ n - bp( G_{n,p})$ becomes less natural as the gap between this
quantity and $ \a( G_{n,p})$ clearly grows.  For example, consider $ p = 1 -q$ where $ q= n^{\gamma}$ and $ -1 \le \gamma < -1/2 $. Let
$Z$ be the number of edges in the largest induced matching in the complement of $ G_{n,p}$. Standard methods show that $ Z =  \Omega( n^{-\gamma} )$ whp (to see this, consider the number of isolated edges in a fixed set of $ 1/q = n^{-\gamma}$ vertices. It is easy to see that the expected number of isolated edges here is linear in $ 1/q$. Furthermore, we can prove concentration around this expected value using the Hoeffding-Azuma inequality).   As the complement of an induced matching with $\ell$ edges 
can be decomposed into $ \ell -1$ bicliques, it follows that 
\[ bp( G_{n,p}) \le n - 2Z+ Z -1 = n - Z -1 \le n -  \Omega( n^{-\gamma} ) \]
whp. On the other hand, the independence number of $ G_{n,p}$ is at most 4 whp in this regime.

\vskip1cm

\noindent
{\bf Acknowledgement:} We thank the anonymous referees for many insightful comments on the original version of this manuscript.

\bibliographystyle{plain}

\end{document}